\documentclass[12pt,reqno]{amsart}

\usepackage{amssymb}
\usepackage{ifthen}
\usepackage{colordvi}

\newcommand{\sbs}{\subset}
\newcommand{\seq}{\subseteq}
\newcommand{\stm}{\setminus}
\newcommand{\est}{\varnothing}

\newcommand{\C}{{\mathbb C}}
\newcommand{\F}{{\mathbb F}}
\newcommand{\K}{{\mathbb K}}
\newcommand{\Q}{{\mathbb Q}}
\newcommand{\R}{{\mathbb R}}

\newcommand{\wh}[1]{{\widehat{#1}}}
\newcommand{\wo}[1]{{\overline{#1}}}

\newcommand{\hG}{\wh{G}}

\newcommand{\hf}{\wh{f}}
\newcommand{\hg}{\wh{g}}
\newcommand{\hw}{\wh{w}}

\newcommand{\prt}{\partial}

\renewcommand{\>}{\rangle}

\DeclareMathOperator{\supp}{supp}

\newcommand{\reft}[1]{\ref{t:#1}}

\newcommand{\refs}[1]{\ref{s:#1}}

\newcommand{\refb}[1]{\cite{b:#1}}
\newcommand{\refe}[1]{\eqref{e:#1}}

\newtheorem{theorem}{Theorem}

\theoremstyle{remark}
\newtheorem{example}{Example}

\newcommand{\showcomments}{no}    

\title[Perfect directions in $\F_p^2$]
  {Point distribution \\ and perfect directions in $\F_p^2$}

\author{Vsevolod F. Lev}
\email{seva@math.haifa.ac.il}
\address{Department of Mathematics, The University of Haifa at Oranim,
  Tivon 36006, Israel}

\begin{document}
\baselineskip = 16pt

\begin{abstract}
Let $p\ge 3$ be a prime, $S\seq\F_p^2$ a nonempty set, and
$w\colon\F_p^2\to\R$ a function with $\supp w=S$. Applying an uncertainty
inequality due to Andr\'as Bir\'o and the present author, we show that there
are at most $\frac12|S|$ directions in $\F_p^2$ such that for every line $l$
in any of these directions, one has
  $$ \sum_{z\in l} w(z) = \frac1p\sum_{z\in\F_p^2} w(z), $$
except if $S$ itself is a line and $w$ is constant on $S$ (in which case all,
but one direction have the property in question). The bound $\frac12|S|$ is
sharp.

As an application, we give a new proof of a result of R\'edei-Megyesi about
the number of directions determined by a set in a finite affine plane.
\end{abstract}

\maketitle

\section{Introduction}\label{s:intro}

Let $p$ be an odd rational prime, and let $\F_p$ denote the $p$-element
field. A \emph{direction} in the affine plane $\F_p^2$ is a pencil of $p$
parallel lines; thus, there are $p+1$ distinct directions.

Given a function $w\colon\F_p^2\to\R$ (which can be thought of as a weight
assignment), we say that a direction is \emph{perfect} with respect to $w$ if
every line in this direction gets its exact share of the total mass of $w$;
that is, for every line $l$ in the direction in question, we have
  $$ \sum_{z\in l} w(z) = \frac1p\,\sum_{z\in\F_p^2} w(z). $$
Write $S:=\supp w$. Choosing a line $l\sbs\F_p^2$ uniformly at random and
considering the variance of the random variable $\sum_{z\in l}w(z)$, it is
easy to show that for all $p+1$ directions to be perfect it is necessary and
sufficient that $w$ be a constant function. Consequently, if all directions
are perfect, then either $S=\F_p^2$, or $S=\est$.%
\ifthenelse{\equal{\showcomments}{yes}}%
{\footnote{\Blue{Assuming that all $p+1$ directions are perfect, let
$w_0:=w-\sum_{x\in\F_p^2} w(x)$. For every line $l$ we then have
 $\sum_{x\in l} w_0(x)=0$ whence, taking the sum over all lines $l\sbs\F_p^2$,
\begin{align*}
  0 &= \sum_l \Big( \sum_{x\in l} w_0(x) \Big)^2 \\
    &= \sum_l \sum_{x,y\in l} w_0(x)w_0(y) \\
    &= \sum_l \sum_{x\in l} w_0^2(x)
          + \sum_{\substack{x,y\in\F_p^2 \\ x\ne y}} w_0(x)w_0(y) \\
    &= p\sum_{x\in\F_p^2} w_0^2(x) + \sum_{x,y\in\F_p^2} w_0(x)w_0(y) \\
    &= p\sum_{x\in\F_p^2} w_0^2(x),
\end{align*}
showing that $w$ is constant.}}
}{} 
In a similar way one can show that a necessary and sufficient condition for
all, but exactly one direction to be perfect is that $w$ is constant on any
line in the unique ``imperfect'' direction; in this case $S$ is a union of
parallel lines, and therefore $|S|\ge p$.

How many perfect directions can there be given that $S$ is small (but
nonempty)? One easily verifies that if $p$ is sufficiently large, then for
$|S|=1$ there cannot be any perfect directions, for $|S|=2$ and $|S|=3$ there
is at most one perfect direction, while for $|S|=4$ there can be two perfect
directions. The goal of this note is to show that, generally, the number of
perfect directions cannot exceed $|S|/2$.
\begin{theorem}\label{t:ud}
Let $p\ge 3$ be a prime. If $S\seq\F_p^2$ is nonempty, then for any function
$w\colon\F_p^2\to\R$ with $\supp w=S$ there are at most $\frac12\,|S|$
perfect directions, unless $S$ is a line and $w$ is constant on $S$ (in which
case there are $p$ perfect directions).
\end{theorem}

A set $S\seq\F_p^2$ is said to \emph{determine} a direction if there is a
line in this direction containing at least two points of $S$. If $|S|=p$ and
$w$ is the indicator function of $S$, then any direction not determined by
$S$ is perfect. Thus, by Theorem~\reft{ud}, if $|S|=p$ and $S$ is not a line,
then there are at most $\frac{p-1}2$ directions not determined by $S$. It
follows that any set $S\seq\F_p^2$ of size $|S|=p$ determines at least
$\frac{p+3}2$ directions, unless $S$ is a line. This is a well-known result
due to R\'edei and Megyesi~\refb{r}, with alternative proofs given by
Lov\'asz and Schrijver~\refb{ls}, and by Dress, Klin, and
Muzichuk~\refb{dkm}. Our Theorem~\reft{ud} thus supplies yet another proof of
this result. In contrast with other proofs, our argument does not rely on the
polynomial method, employing Fourier analysis instead.

We refer the reader to~\refb{g} for a historical account and summary of
related results.

The following examples show that the estimate of Theorem~\reft{ud} is, in a
sense, best possible.

\begin{example}
The special orthogonal group $\mathrm{SO}(2,p)$ is cyclic of order
$p-(-1/p)$, where $(\cdot/p)$ is the Legendre symbol. Assuming that $2n$ is
an even integer dividing $p-(-1/p)$, let $H\le\mathrm{SO}(2,p)$ be the
subgroup of order $|H|=2n$, and let $H_0<H$ be the subgroup of $H$ of order
$|H_0|=n$. Fix arbitrarily a nonzero point $z\in\F_p^2$, define $S$ to be the
orbit of $z$ under the action of $H$, and for $x\in S$ let $w(x)=1$ if $x$
actually belongs to the orbit of $z$ under the action of $H_0$, and $w(x)=-1$
otherwise. We leave it to the reader to verify that there are $n=\frac12|S|$
directions determined by the pairs $(x,y)\in S\times S$ with $w(x)\ne w(y)$,
and that all these directions are perfect.
\end{example}

The next example originates, essentially, from Lov\'asz-Schrijver~\refb{ls}.
\begin{example}
Let $S$ be the graph of the function $z\mapsto z^{\frac{p+1}2}$; that is,
$S=\{(z,z^{\frac{p+1}2})\colon z\in\F_p\}$. Then $S$ determines $\frac{p+3}2$
directions, and since $|S|=p$, the $\frac{p-1}2=\lfloor|S|/2\rfloor$
undetermined directions are perfect with respect to the indicator function of
$S$.
\ifthenelse{\equal{\showcomments}{yes}}{%
\Blue{(Indeed, as shown in~\refb{ls}, up to affine transformations, $S$ is
the unique set of size $|S|=p$ determining exactly $\frac{p+3}2$
directions.)}
}{}
\end{example}

\begin{example}
If $l_1,l_2\sbs\F_p^2$ are nonparallel lines, and $w$ is the difference of
the indicator functions of these lines, then
 $S=(l_1\cup l_2)\stm(l_1\cap l_2)$, $|S|=2(p-1)$, and there are $p-1=|S|/2$
perfect directions. Similarly, if $S$ is a union of two \emph{parallel}
lines, and $w$ is constant and nonzero on each of these lines, then $|S|=2p$
and there are $p=|S|/2$ perfect directions.
\end{example}

We prove Theorem~\reft{ud} in the next section, and discuss related open
problems in the concluding Section~\refs{conclusion}.

\section{The proof of Theorem~\reft{ud}}\label{s:proof}

We begin with setting up the notation and recalling basic facts and
properties of the Fourier transform on finite abelian groups.

For a subfield $\K$ of the field $\C$ and a finite, nonempty set $G$, by
$L_\K(G)$ we denote the space of all functions from $G$ to $\K$ with the
inner product defined by
  $$ \<f,g\> := \frac1{|G|} \sum_{z\in G} f(z)\,\wo{g(z)},\quad f,g\in
                                                                  L_\K(G), $$
the overline denoting the complex conjugation.

Suppose that $G$ is a finite abelian group. \emph{Dual} to $G$ is the group
of all homomorphisms from $G$ to $\C^\times$. The dual group is denoted
$\hG$, its elements are called \emph{characters}, the identity element of
$\hG$ is the \emph{principal} character. The Fourier transform of a function
$f\in L_\K(G)$ is the function $\hf\in L_\C(\hG)$ defined by
  $$ \hf(\chi) := \<f,\chi\>, \quad \chi\in \hG. $$
The function $f\in L_\K(G)$ is constant if and only if its Fourier transform
is zero or supported on the principal character.

For a subgroup $H\le G$, the set of all characters $\chi\in\hG$ containing
$H$ in their kernel is a subgroup of $\hG$, denoted $H^\perp$; thus,
  $$ H^\perp
      = \{ \chi\in\hG\colon \chi(h)=1\ \text{for any}\ h\in H \} \le \hG. $$
If $H\le G$ is nonzero and proper, then so is $H^\perp\le\hG$. Writing $1_H$
and $1_{H^\perp}$ for the indicator functions of $H$ and $H^\perp$,
respectively, we have $\wh{1_H}=(|H|/|G|)\cdot1_{H^\perp}$.

For a function $g\in L_\K(G)$ and an element $z\in G$, let $g_z\in L_\K(G)$
be defined by
   $$ g_z(x) := \wo{g(z-x)},\quad x\in G. $$
The \emph{convolution} of functions $f,g\in L_\K(G)$ is the function
  $$ f\ast g \colon z \mapsto \<f,g_z\>, \quad z\in G. $$
The Fourier transform of a convolution is the product of Fourier transforms:
  $$ \wh{f\ast g} = \hf\cdot\hg,\quad f,g\in L_\K(G). $$

Our argument relies on the following uncertainty inequality for the
rational-valued functions on the finite affine plane.
\begin{theorem}[{\cite[Theorem~1]{b:bl}}]\label{t:bl}
For any prime $p\ge 3$ and any function $f\in L_\Q(\F_p^2)$, either
  $$ \frac12\,|\supp f| + \frac1{p-1}\,|\supp\hf| \ge p+1, $$
or there is a direction in $\F_p^2$ such that $f$ is constant on every line
in this direction.
\end{theorem}

We now turn to the proof of Theorem~\reft{ud}.

If $w$ is a constant function, then $S=\F_p^2$ and the assertion is
immediate; assume thus that $w$ is not constant. The case $p=3$ is easy to
verify, and we further assume that $p\ge 5$.

By the Dirichlet simultaneous approximation theorem, there exist arbitrarily
large integers $Q$, along with the corresponding integer-valued functions
$w_Q$ on $\F_p^2$, such that
  $$ \left\| w - \frac{w_Q}Q\right\|_\infty < \frac1{2pQ}. $$
As a result, if $Q$ is sufficiently large, then $\supp w_Q =\supp w$, and for
$x,y\in S$ we have $w(x)=w(y)$ if and only if $w_Q(x)=w_Q(y)$; also, a
direction is perfect with respect to $w$ if and only if it is perfect with
respect to $w_Q$. Consequently, passing from $w$ to $w_Q$, we can ensure
that, in addition to being nonconstant, $w$ is also integer-valued.

To every direction in $\F_p^2$ there corresponds a nonzero, proper subgroup
$H<\F_p^2$; specifically, the subgroup represented by the line through the
origin in the corresponding direction. As an immediate corollary from the
definitions, the direction corresponding to a subgroup $H<\F_p^2$ is perfect
if and only if the convolution $w\ast 1_H$ is a constant function; that is,
the product $\widehat w\cdot 1_{H^\perp}$ vanishes at every nonprinciple
character; in other words, $\widehat w$ vanishes on every character from
$H^\perp$ with the possible exception of the principle character.
\ifthenelse{\equal{\showcomments}{yes}}%
{\Blue{As an immediate corollary, we have yet another explanation for the
claims at the beginning of the Introduction:
\begin{itemize}
\item[i)]   if $w$ is not constant, then there is at least one imperfect
    direction;
\item[ii)]  the direction corresponding to a subgroup $H<\F_p^2$ is the
    unique imperfect direction if and only if $w$ is $H$-periodic; that
    is, $w(g+h)=w(g)$ for any $g\in\F_p^2$ and $h\in H$.
\end{itemize}}
}{} 

Denote the number of perfect directions by $N$, so that the number of
imperfect directions is $p+1-N$. The group $\wh{\F_p^2}\cong\F_p^2$ is a
union of its $p+1$ nonzero, proper subgroups, with every nonprincipal
character $\chi\in\wh{\F_p^2}$ lying in exactly one subgroup, and the
principal character lying in all subgroups. Therefore, since $\hw$ vanishes
on the subgroups corresponding to the perfect directions, we have
\begin{equation}\label{e:theequation}
  |\supp\wh{w}| \le (p-1)(p+1-N) + 1.
\end{equation}

On the other hand, applying Theorem~\reft{bl} to the function $w$, we
conclude that either
\begin{equation}\label{e:lastone}
  \frac12\,|\supp w| + \frac1{p-1}\,|\supp\widehat w| \ge p+1,
\end{equation}
or there is a direction $\prt$ such that $w$ is constant on every line in
this direction. In the former case, combining~\refe{lastone}
and~\refe{theequation}, and recalling that $p\ge 5$, we get
  $$ p+1 \le \frac12\,|S| + \left( (p+1-N) + \frac1{p-1} \right)
                < \frac12\,(|S|+1) + p+1-N, $$
implying $N\le\frac12\,|S|$. In the latter case, denoting by $k$ the number
of lines in the direction $\prt$ on which $w$ is nonzero, we have $|S|=kp$,
while $N=p$ (all directions except $\prt$ are perfect). Consequently,
$N\le\frac12\,|S|$, unless $k=1$, meaning that there is a line on which $w$
is constant and nonzero, and outside of which $w$ vanishes.

This completes the proof of Theorem~\reft{ud}.

\section{Open problems: restricting the weights}\label{s:conclusion}

Suppose that $w\in L_\C(\F_p^2)$ is not constant, and let $S=\supp w$. If
$|S|<p$, then in every direction there is a line disjoint from $S$; hence,
for perfect directions to exist, the average value of $w$ on $\F_p^2$ must be
zero. This suggests the following problem: how many perfect directions can
there be for a function $w$ with a small support given that the average of
$w$ is \emph{nonzero}?

As we have just saw, one needs $|S|\ge p$ in order to have any perfect
directions at all.

If $|S|=p$, 
then for any direction determined by $S$ there is a line in this direction
disjoint from $S$; therefore, none of the directions determined by $S$ is
perfect. On the other hand, if $w$ is constant on $S$, then any direction not
determined by $S$ is perfect. It follows that the largest possible number of
perfect directions is equal to the largest possible number of undetermined
directions, which is $p+1$ less the smallest possible number of determined
directions. Apart from the trivial case where $S$ is a line, the smallest
possible number of determined directions is $(p+3)/2$ by the of
R\'edei-Megyesi result; thus, the largest possible number of perfect
directions (for $S$ not being a line, $|S|=p$, and the average of $w$
nonzero) is $(p+1)-\frac12(p+3)=\frac12(p-1)$.

For $|S|=p+1$, one perfect direction is very easy to arrange, and a simple
combinatorial argument shows that there cannot be two or more perfect
directions. Notice that this contrasts sharply the situation where $|S|=p$.

For $|S|=p+2$ one can have two perfect directions (set $w(z)=1/2$ for
$z\in\{(0,0),(0,1),(1,0),(1,1)\}$, and $w(z)=1$ for $z=(x,x)$ with
$x\in[2,p-1]$); however, it is not clear to us whether there can be three or
more perfect directions, nor what happens for $|S|\ge p+3$.

Replacing the nonzero average assumption with the stronger assumption that
$w$ attains real \emph{nonnegative} values seems to result in an equally
interesting problem.

\bigskip

\end{document}